\documentstyle{amsppt}
\NoBlackBoxes

\define\cc{\Bbb C}
\define\cn{\Bbb C^n}
\define\cm{\Bbb C^m}
\define\cd{\Bbb C^2}

\define\er{\Bbb R}
\define\ro{\varrho}
\define\wl{\Cal L_0}
\define\ord{\operatorname{ord}}
\define\OO{\Cal O}
\define\la{\lambda}
\define\vf{\varphi}

\define\oo#1{\overset {{}_o}\to {#1}}
\define\zo{\overset {{}_o}\to z}

\document
\topmatter 
\footnote"{}"{1991 {\it Mathematics Subject Classification}: 32S05.} 
\footnote"{}"{{\it Key words and phrases}: holomorphic mapping, \L ojasiewicz
exponent.} 
\footnotetext"{}"{This research was partially supported by KBN Grant No. 2
P03A 050 10.}
\endtopmatter

\vglue 1.5cm
\centerline{\bf A set on which the local \L ojasiewicz exponent}

\centerline{\bf is attained}

\vskip.5cm
\font\tw=cmcsc10
\centerline{\tw {\rm by} Jacek Ch\c adzy\'nski and Tadeusz Krasi\'nski {\rm (\L
\'od\'z)}} 


\font\bbb=cmr8
\font\ccc=cmbx8

\vskip.7cm
\vbox{{\ccc Abstract. }\bbb Let $U$ be a neighbourhood of $0\in \cn$. We show
that for a holomorphic mapping $F=(f_1,...,f_m): U\to\cm$, $F(0)=0$, the \L
ojasiewicz exponent $\wl(F)$ is attained on the set $\{z\in
U:\,f_1(z)\cdot \ldots \cdot f_m(z)=0\}$.}

\vskip.5cm 
{\bf 1. Introduction.} In the paper [CK$_2$] the autors showed that for a
polynomial mapping $F=(f_1,...,f_m):\cn\to\cm$, $n\ge 2$, the \L ojasiewicz
exponent $\Cal L_\infty(F)$ of $F$ at infinity is attained on the set
$\{z\in\cn:f_1(z)\cdot\ldots\cdot f_m(z)=0\}$. The purpose of this paper is
to prove an analogous result for the \L ojasiewicz exponent $\wl(F)$, where
$F:U\to \cm$ is a holomorphic mapping, $F(0)=0$ and $U$ is a neighbourhood
of $0\in\cn$ (Thm 1). From this result we easily obtain a strict formula for
$\wl(F)$ in the case $n=2$ and $m\ge 2$ in terms of multiplicities of some
mappings from $U$ into $\cd$ defined by components of $F$ (Thm 2). It is a
generalization of Main Theorem from [CK$_1$]. Proof of this theorem has been
simplified by A. P\l oski in [P]. His proof has been an inspiration to write
this paper.

Theorem 1 is an important tool for investigation of the \L ojasiewicz
exponent for analytic curves having an isolated intersection point at
$0\in\cm$. Using it, we shall give in the next paper [CK$_3$], an effective
formula for the \L ojasiewicz exponent for such curves in terms of their
parametrizations. 

{\bf 2. The \L ojasiewicz exponent.} Let $U\subset \cn$, $n\ge 2$, be a
neighbourhood of the origin, $F:U\to \cm$ a holomorphic mapping, $S\subset U$
an analytic set in $U$. Assume that $0\in\cn$ is an accumulation point of
$S$. Put
$$
N(F|S):=\{\nu \in\er_+:\exists A>0,\exists B>0, \forall z\in S,
|z|<B\Rightarrow A|z|^\nu\le |F(z)|\}.
$$
Here $|\cdot |$ means the polycylindric norm. If $S=U$ we write $N(F)$
instead of $N(F|U)$.

By the \L ojasiewicz exponent of the mapping $F|S$ at $0$ we mean
$$
\wl(F|S):=\inf N(F|S).
$$
Analogously $\wl(F):=\inf N(F)$.

It can be shown (cf [LT], \S\S 5,6)

\proclaim{\indent{\tw Proposition 1}} If $F|S$ has an isolated zero at
$0\in\cn$, then $\wl(F|S)\in N(F|S)\cap \Bbb Q$. Moreover, there exists an
analytic complex curve $\vf:\{t\in\cc:|t|<r\}\to S\}$ such that $\vf(0)=0$ and
$$
|F\circ \vf(t)|\sim |\vf(t)|^{\wl(F|S)}\qquad\text{for }t\to 0.
\tag 1
$$
\endproclaim

From the above proposition we easily get

\proclaim{\indent{\tw Corollary 1}} $\wl(F|S)<\infty$ if and only if $F|S$
has an isolated zero of $0\in \cn$.
\endproclaim

{\bf 3. The main result.} Now, we shall give the main result of this paper.

\proclaim{\indent{\tw Theorem 1}} Let $U\subset \cn$, $n\ge 2$, be a
neighbourhood of the origin, $F=(f_1,...,f_m):U\to \cm$ a holomorphic
mapping, and $F(0)=0$. Define $S:=\{z\in U:f_1(z)\cdot\ldots\cdot
f_m(z)=0\}$. Then
$$
\wl(F)=\wl(F|S).
\tag 2
$$
\endproclaim

The proof will be given in section 4.

Immediately from Theorem 1 we obtain an effective formula for the \L
ojasiewicz exponent in the case $n=2$, $m\ge 2$, generalizing an earlier
authors' result ([CK$_1$], Main Thm).

Let us begin with notations. Let $U$ be a neighbourhood of $0\in \cd$. Then:
$\mu(f,g)$ means the intersection multiplicity of a holomorphic mapping
$(f,g):U\to \cd$; $\hat h$ means the germ of a holomorphic function $h:U\to
\cc$ in the ring $\OO^2$ of germs of holomorphic functions at $0\in\cd$;
$\ord h$ means the order of $h$ at $0\in\cd$.

\proclaim{\indent{\tw Theorem 2}} Let $U\subset \cd$ be a neighbourhood of
the origin, $F=(f_1,...,f_m):U\to \cm$ a holomorphic mapping, $F(0)=0$. Put
$f:=f_1\cdot\ldots\cdot f_m$. If $\hat f\ne 0$ and $\hat f=\hat
h_1\cdot\ldots\cdot \hat h_r$ is a factorization of $\hat f$ into irreducible
germs in $\OO^2$, then
$$ 
\wl(F)=\max_{i=1}^r \frac 1{\ord h_i} \min_{j=1}^m \mu(h_i,f_j).
\tag 3
$$
\endproclaim

\demo{Proof} Since our considerations are local, we may assume that $h_i$ are
holomorphic in $U$ and $f=h_1\cdot\ldots\cdot h_r$ in $U$. Let $S:=\{z\in
U:f(z)=0\}$ and $\Gamma_i:=\{z\in U:h_i(z)=0\}$. Hence
$$
S=\Gamma_1\cup\ldots\cup \Gamma_r.
\tag 4
$$

Define $\la_i:=(1\slash \ord h_i)\min_{j=1}^m\mu(h_i,f_j)$. If
$\la_i=+\infty$ for some $i$, then (3) holds. So, assume that $\la_i\ne
+\infty$, $i=1,...,r$. Then for every $i$ we have
$$
|F(z)|\sim |z|^{\la_i}\qquad\text{for }|z|\to 0\text{ and }z\in\Gamma_i.
$$
Hence and from (4), $\wl(F|S)=\max_{i=1}^r \la_i$. This and Theorem 1 give
(3), which ends the proof.
\enddemo

{\bf 4. Proof of the main theorem.} The proof is given in two steps. In
the first one we give the proof under some additional assumptions,
whereas in the second one, we show that these assumptions do not
restrict our considerations.

First we fix some notations. For $z=(z_1,...,z_n)\in \cn$, $n\ge 2$, and for
every $i\in\{1,...,n\}$ we put $z_i':=(z_1,...,z_{i-1},z_{i+1},...,z_n)$.
Additionally, we define $f:=f_1\cdot\ldots\cdot f_m$.

{\it Step 1.} We assume that

(i) $(F|S)^{-1}(0)=\{0\}$,

(ii) $\ord f<\infty$,

(iii) for every $i\in\{1,...,n\}$, $f$ is $(\ord f)$-regular with respect to
$z_i$. 

Obviously $N(F)\subset N(F|S)$. To show (2) it suffices to prove
$$
N(F|S)\subset N(F).
\tag 5
$$

It follows from (i) and Corollary 1 that $N(F|S)$ is not empty. Take an
arbitrary $\nu\in N(F|S)$. Then there exist $A,\,B>0$ such that
$$
|F(\zeta)|\ge A|\zeta |^\nu\qquad\text{for }\zeta \in S,\; |\zeta |<B.
\tag 6
$$

From the assumptions of the theorem we have $\ord f>0$ and $\ord f_j>0$,
$j=1,...,m$. Then from (ii), (iii) we easily get $0<\ord f_j<\infty$ and 
$$
\text{for every }i,\,j,\; f_j\text{ is }(\ord f_j\text{)-regular with respect
to }z_i.
\tag 7
$$
The Weierstrass Preparation Theorem and (iii) imply that for every
$i\in\{1,...,n\}$ there exists a Weierstrass polynomial with respect to $z_i$
of degree $\ord f$, associated with $f$. Denote it by $p^{(i)}$. Analogously
from (7) for every $i$, $j$ there exists a Weierstrass  polynomial with
respect to $z_i$ of degree $\ord f_j$, associated with $f_j$. Denote it by
$p_j^{(i)}$. Then there exist $C,\,D>0$ and a polycylinder $\{z\in \cn:
|z|<r\}\subset U$ such that
$$
C|p^{(i)}(z)|\le |f(z)|\le D|p^{(i)}(z)|\qquad\text{for }|z|<r,\;i=1,...,n,
\tag 8
$$
and
$$
C|p_j^{(i)}(z)|\le |f_j(z)|\le D|p_j^{(i)}(z)|\qquad\text{for
}|z|<r,\;i=1,...,n,\,j=1,...,m. 
\tag 9
$$
Put $P^{(i)}:=(p_1^{(i)},...,p_m^{(i)})$. Then from (9) we get
$$
C|P^{(i)}(z)|\le |F(z)|\le D|P^{(i)}(z)|\qquad\text{for }|z|<r,\,i=1,...,n.
\tag 10
$$
Clearly, for every $i\in\{1,...,n\}$
$$
p^{(i)}=p_1^{(i)}\cdot\ldots\cdot p_m^{(i)}.
\tag 11
$$
Put $B_1:=\min (B,r)$. Then from (8) we have for every $i\in\{1,...,n\}$
$$
S\cap \{z\in\cn:|z|<B_1\}=\{z\in \cn:|z|<B_1,\,p^{(i)}(z)=0\}.
$$
Hence and from the Theorem on Continuity of Roots, applied to $p^{(i)}$, we
get that there exists $\ro$, $0<\ro\le B_1$, such that for every
$i\in\{1,...,n\}$ we have
$$
\{z\in S:|z_i'|<\ro,\,|z|<B_1\}=\{z\in \cn:|z_i'|<\ro,\,p^{(i)}(z)=0\}.
\tag 12
$$

Put $d:=\max_{j=1}^m \ord f_j$, $A_2:=2^{-d}A(C\slash D)$, $B_2:=\min
(\ro,B_1)$. Take an arbitrary $\zo\in \cn$, $|\zo|<B_2$. There exists $i\in
\{1,...,n\}$ such that $|\zo|=|\zo_i'|$. Define
$\vf_j(t):=p_j^{(i)}(\zo_1,...,\zo_{i-1},t,\zo_{i+1},...,\zo_n)$,
$\Phi:=(\vf_1,...,\vf_m)$. Clearly, $\Phi:\cc\to\cm$ is a polynomial mapping,
$\deg \Phi:=\max_{j=1}^m\deg \vf_j=d$ and
$\Phi(t)=P^{(i)}(\zo_1,...,\zo_{i-1},t,\zo_{i+1},...,\zo_n)$. Then by Lemma
2 in [CK$_2$] and (11) we have
$$
|\Phi(\zo_i)|\ge 2^{-d}\min_{\tau\in T}|\Phi(\tau)|,
$$
where $T:=\{t\in \cc:p^{(i)}(\zo_1,...,\zo_{i-1},t,\zo_{i+1},...,\zo_n)=0\}$.
Hence and from (10) we get
$$
\aligned
|F(\zo)|&\ge C|P^{(i)}(\zo)|=C|\Phi(\zo_i)|\ge C2^{-d}|\Phi(\tau _0)|\\
&=C2^{-d}|P^{(i)}(\oo{\zeta})|\ge (C\slash D)2^{-d}|F(\oo{\zeta})|
\endaligned
\tag 13
$$
for some $\oo{\zeta}=(\zo_1,...,\zo_{i-1},\tau_0,\zo_{i+1},...,\zo_n)$ such
that $p^{(i)}(\oo{\zeta})=0$. Since $|\oo{\zeta}_i'|=|\zo_i'|=|\zo|<B_2\le
\ro$, then from (12) we have $\oo{\zeta}\in S$ and $|\oo{\zeta}|<B_1$. In
consequence, from (6) and (13) we get
$$
|F(\zo)|\ge (C\slash D)2^{-d}A|\oo{\zeta}|^\nu\ge
A_2|\oo{\zeta}_i'|^\nu=A_2|\zo|^\nu. 
$$
Since $\zo$ is arbitrary we obtain $\nu\in N(F)$. This ends the proof of the
theorem under assumptions (i), (ii), (iii).

{\it Step 2.} Now, we shall prove the theorem in the remaining
cases. If (i) does not hold, then by Corollary 1 we have
$\wl(F|S)=\wl(F)=\infty$, that is, (2) is satisfied. If (ii) does not hold,
then $S=U$ and (2) is obvious. So, it suffices to consider the case when (i),
(ii) hold but (iii) does not. Since $N(F)$ and $N(F|S)$ are invariant with respect
to linear automorphisms of $\cn$, then $\wl(F)$ and $\wl(F|S)$ are also
invariant of such automorphisms. Since the considered case can be reduced to
step 1 by a linear automorphism of $\cn$ then (2) also holds in this case.

This ends the proof of the theorem.

\font\bbb=cmr8

\vskip.3cm
\noindent{\baselineskip=10pt
\bbb Faculty of Mathematics\newline
University of \L \'od\'z\newline
S. Banacha 22\newline
90-238 \L \'od\'z, Poland

\noindent
E-mail: jachadzy\@imul.uni.lodz.pl

\noindent\hskip 11mm krasinsk\@krysia.uni.lodz.pl}
\enddocument